\documentclass{amsart}
\usepackage{amssymb,colonequals,enumerate}
\usepackage{tikz-cd}

\usepackage{hyperref}
\hypersetup{%
  bookmarksnumbered=true,%
  colorlinks=true,%
  linkcolor=blue,%
  citecolor=blue,%
  filecolor=blue,%
  menucolor=blue,%
  urlcolor=blue,%
  bookmarksopen=true,%
  bookmarksdepth=2,%
  pageanchor=true}

\numberwithin{equation}{section}

\swapnumbers

\theoremstyle{plain}
\newtheorem{theorem}[equation]{Theorem}
\newtheorem{proposition}[equation]{Proposition}
\newtheorem{lemma}[equation]{Lemma} 
 
\newtheorem{corollary}[equation]{Corollary}

\theoremstyle{definition}
\newtheorem{definition}[equation]{Definition}
\newtheorem{example}[equation]{Example}

\newtheorem{chunk}[equation]{}

\theoremstyle{remark}
 
\newtheorem{question}[equation]{Question}
\newtheorem*{ack}{Acknowledgements}

\newcommand{\cat}[1]{\mathsf{#1}}

\newcommand{\depth}{\operatorname{depth}}
\renewcommand{\dim}{\operatorname{dim}}

\newcommand{\dcat}[2][{}]{{\cat D}_{#1}(#2)}

\newcommand{\Ann}{\operatorname{Ann}}
\newcommand{\edim}{\operatorname{edim}}
\newcommand{\End}{\operatorname{End}}

\newcommand{\fdim}{\operatorname{flat\,dim}}
\newcommand{\fm}{\mathfrak{m}} 
\newcommand{\hh}{\operatorname{H}}

\newcommand{\id}{\operatorname{id}}
\newcommand{\lotimes}{\otimes^{\operatorname{L}}}
\newcommand{\op}[1]{{#1}^{\operatorname{op}}}
\newcommand{\pdim}{\operatorname{proj\,dim}}
\newcommand{\perf}[2][{}]{\operatorname{Perf}_{#1}(#2)}
\newcommand{\PGL}{\mathrm{PGL}}
\newcommand{\rank}{\operatorname{rank}}
\newcommand{\RHom}{\operatorname{RHom}}

\newcommand{\Tor}{\operatorname{Tor}}

\newcommand{\vf}{\varphi}

\newcommand{\bs}{\boldsymbol}

\newcommand{\mcO}{\mathcal{O}}

\newcommand{\mbb}{\mathbb}

\newcommand{\rhobar}{\overline{\rho}}

\begin{document}

\title[Freeness criterion]{A freeness criterion for complexes\\ with derived actions}

\author[Brochard]{Sylvain Brochard}
\address{IMAG, University of Montpellier, CNRS, Montpellier, France}
\email{sylvain.brochard@umontpellier.fr}

\author[Iyengar]{Srikanth B.~Iyengar}
\address{Department of Mathematics,
University of Utah, Salt Lake City, UT 84112, U.S.A.}
\email{srikanth.b.iyengar@utah.edu}

\author[Khare]{Chandrashekhar  B. Khare}
\address{Department of Mathematics,
University of California, Los Angeles, CA 90095, U.S.A.}
\email{shekhar@math.ucla.edu}

\date{\today}


\begin{abstract} 
Inspired by the patching method of Calegari and Geraghty, and a conjecture of de Smit that has been proved by the first author, we present a conjectural freeness criterion without patching for complexes over commutative noetherian local rings with derived actions, and verify it in several cases. 
\end{abstract}

\maketitle
   
\section{Introduction} 

This work is motivated by a patching method due to  Calegari and Geraghty \cite{Calegari/Geraghty:2018} that  extends the results in  \cite{Taylor/Wiles:1995} to situations where one patches complexes rather than modules, as we explain below. We use  the reformulation of the method of \cite{Calegari/Geraghty:2018} as given in \cite{Khare/Thorne:2017}.

Let $\mcO$ be a complete DVR and $F=0\to F_p\to \cdots \to F_0\to 0$ a complex of  free $\mcO$-modules of finite rank with $\hh_0(F)\ne 0$; $p $ is called the \emph{defect}. Let $B$ be a complete noetherian local $\mcO$-algebra equipped with a morphism of $\mcO$-algebras $B \to \End_{\cat D(\mcO)}(F)$, where $\cat D(\mcO)$ is the derived category of $\mcO$-modules. One wants to find a criterion for proving that $\hh_0(F)$ is a free $B$-module. 

The  following ``patching data'' in  \cite{Khare/Thorne:2017}  gives such a criterion.
For each integer  $n\ge 1$  suppose  one is given:
 \begin{enumerate}[\rm(i)]
 \item
 Complete noetherian local $\mcO$-algebras $A_n,B_n$  with  $\edim A_n-\edim B_n=p $, where $p $ is as above; in particular, it is independent of $n$.
 \item
 Local maps $\phi_n\colon A_n \to B_n$ and   complexes $F_n$ of  length $\leq p $  consisting  of free $A_n$-modules of  rank bounded independently of $n$.
 \item A morphism $B_n \to \End_{\cat D(A_n)}(F_n)$. 
 \end{enumerate}
 It is also assumed that $A_n$ somehow ``grows’’ with $n$; for the precise condition we refer the reader to   \cite[Proposition 3.1, Theorem 6.29]{Khare/Thorne:2017}.
 
Patching this data, which relies on a compactness argument  and uses the pigeonhole principle,  leads to $\mcO$-algebras $A_\infty$ and $B_\infty$, with $A_\infty$ a regular local ring,  and a complex $F_\infty$ of finite free  $A_\infty$-modules of length $\le p $ with the property that $F_\infty$ is exact: $\hh_i(F_\infty)=0$ for $i\ge 1$. Using the Auslander-Buchsbaum formula, it is easy to prove that $\hh_0(F_\infty)$ is a free $B_\infty$-module. One recovers $B$ and $F$  from $B_\infty$ and $F_\infty$ by going modulo the kernel of a section $A_\infty \to \mcO$ of the structure map $\mcO \to A_\infty$.  This implies that $\hh_0(F)$ is a free $B$-module.

The  question below seeks to obviate the need to patch in \cite{Calegari/Geraghty:2018}, by asking for the patching data only for $n=1$,  just as de Smit's conjecture, proved in~\cite{Brochard:2017}, obviated the need to patch in the original argument of Taylor and Wiles. 

\begin{question}
\label{qu:Shekhar}
Let $A$ be a commutative noetherian local ring and 
\[
F\colon 0\longrightarrow F_p\longrightarrow \cdots \longrightarrow F_0\longrightarrow 0
\]
a finite free $A$-complex with  $\hh_0(F)\ne 0$. If the canonical map $A\to \End_{\cat D(A)}(F)$ factors through a commutative noetherian local ring $B$ such that
\[
p \le  \edim A - \edim B \,,
\]
is then the $B$-module $\hh_{0}(F)$ free?
\end{question}

As noted above, we know since the work of Calegari and Geraghty that Question~\ref{qu:Shekhar} has a positive answer if the ring $A$ is regular; see~\cite[Section 6]{Calegari/Geraghty:2018}. The results of \cite{Brochard:2017, Brochard/Iyengar/Khare:2023a} settle the case where $\hh_i(F)=0$ for $i\ne 0$, that is to say, when $F$ is quasi-isomorphic to a $B$-module; see also Theorem~\ref{th:main1} below. We are unable to verify Question~\ref{qu:Shekhar} in full generality. The goal of this manuscript is more modest: to explore the viability of the question by testing the tightness of the hypotheses and establishing an affirmative answer to it in some cases. 

For a start, we remark that a derived $B$-action on $F$ induces a $B$-action on $\hh_*(F)$, the homology of $F$. However, for the conclusion of Question~\ref{qu:Shekhar} to hold it does not suffice that $\hh_*(F)$ is a $B$-module; see Example~\ref{ex:H-action}, but also Theorem~\ref{th:regular_case}.  

On the other hand, imposing the more stringent hypotheses that $F$ is a complex of $B$-modules---equivalently, that there is a homomorphism  $B\to \End_A(F)$ of $A$-algebras, where $\End_A(F)$ is the endomorphism differential graded algebra---turns out to be too restrictive; this is the import of the result below. Thus in summary  a derived $B$-action on $F$  has the Goldilocks-like quality of  being  neither  too weak nor too strong, making us hopeful  that the question  has an  affirmative answer without it  leading to conclusions  that are untenable in  number theoretic situations of the kind that arise in patching complexes rather than modules.

\begin{theorem}
 \label{th:main1}
 In the context of Question~\ref{qu:Shekhar} assume furthermore that  $F$ is quasi-isomorphic in $\dcat A$ to a complex of $B$-modules. Then $\hh_i(F)=0$ for $i\ne 0$, the $B$-module $\hh_0(F)$ is free, and the map $A\to B$ is an exceptional complete intersection.
\end{theorem}

This result, which is proved in Section~\ref{section:B_modules}, is a generalization, with a new proof, of the main result of \cite{Brochard/Iyengar/Khare:2023a}, that dealt with the case where it was assumed a priori that $\hh_i(F)=0$ for $i\ne 0$, that is to say that $F$ is quasi-isomorphic to a $B$-module. The real import of the theorem  from the perspective of number theoretic applications is that tightening the hypotheses on the complex $F$ in Question~\ref{qu:Shekhar} is not reasonable. Indeed, in Section~\ref{section:nt} we explain that  as a consequence the theorem,  the ``patching data’’  cannot  be strengthened  to   asking for morphisms $B_n \to \End_{A_n}(F_n)$ in most  situations of interest.  In a different direction, one may ask if  assuming  $F$ has Poincare duality--an assumption which is fulfilled in the number theoretic situations  (see  Section~\ref{section:nt}) in which the  complexes  compute the homology of arithmetic manifolds--helps in answering the question.

The result below is proved in Section~\ref{section:ci_case}. 

\begin{theorem}
 \label{th:main2}
 Question~\ref{qu:Shekhar} has a positive answer when $B/\fm_AB$ is a complete intersection, and in particular when the map $A\to B$ is surjective.
\end{theorem}

Different ideas are needed to prove it than those that go into proving Theorem~\ref{th:main1}. For instance, the first step in the proof of \emph{op.~cit.} is a reduction to the case where the map $A\to B$ is surjective; this is one place where the hypothesis that $F$ is equivalent to a $B$-complex is crucial. No such reduction seems possible in proving the theorem above.  The argument relies instead on a freeness criterion for balanced modules by the first author~\cite{Brochard:2023}, and on a characterization of graded modules over a certain kind of non-commutative Weyl algebra given in Appendix~\ref{se:appendix}.  Moreover, unlike in Theorem~\ref{th:main1}, there are complexes $F$ that satisfy the hypotheses of the preceding theorem, having homology in multiple degrees; moreover, the map $A\to B$ need not be a complete intersection; see Example~\ref{ex:koszul}. We see this as good evidence that the question above might have an affirmative answer.



\begin{ack}
    This work is partly supported by National Science Foundation grants DMS-200985 (SBI) and  DMS-2200390 (CBK).   We would  also like to thank  Najmuddin Fakhruddin and Jack Thorne for helpful discussions.

\end{ack}

\section{Preliminaries}
\label{section:preliminaries}

In this work by a ring we mean a commutative ring, unless stated otherwise. Throughout the paper $A$ and $B$ denote noetherian local rings. We write $\dcat A$ for the derived category of $A$-complexes. For any $M$ in $\dcat A$ set
\[
\inf\hh_*(M) =\inf\{n\mid \hh_n(M)\ne 0\}\quad\text{and}\quad \sup\hh_*(M) =\sup\{n\mid \hh_n(M)\ne 0\}\,.
\]
We begin by recalling some notions and facts concerning homological invariants of complexes over local rings; see  Foxby~\cite{Foxby:1979a}.

\subsection*{Homological invariants}
Let $A$ be a noetherian local ring with maximal ideal $\fm_A$ and residue field $k_A$. The \emph{depth} of an $A$-complex $M$ is 
\[
\depth_AM \colonequals \inf\{n\mid \mathrm{Ext}^n_A(k_A,M)\ne0\}\,,
\]
and the \emph{dimension} of an  $A$-complex $M$ is 
\[
\dim_AM \colonequals \sup\{\dim_A \hh_n(M) - n \mid n\in \mathbb{Z}\}\,.
\]
The \emph{Cohen-Macaulay defect} of $M$ is $\dim_AM - \depth_AM$.

Let $M$ be an $A$-complex such that the $A$-module $\hh_*(M)$ is finitely generated; that is to say, the $A$-module $\hh_i(M)$ is finitely generated for each $i$, and equal to zero when $|i|\gg 0$. The \emph{projective dimension} of $M$ is 
\[
\pdim_AM \colonequals \sup\{n\mid \mathrm{Ext}^n_A(M,k_A)\ne 0\} = \sup\{n\mid \Tor^A_n(k_A,M)\ne 0\}\,.
\]
When the projective dimension of $M$ is finite  we say it is \emph{perfect}; this is equivalent to the condition that $M$ is isomorphic in $\dcat A$ to a finite free complex.
For this reason, we often use the letter $F$ to denote a perfect complex. Note that if a finite free complex $F$ is \emph{minimal}, that is to say, if $\partial(F)\subseteq \fm_A F$, then its projective dimension is the largest $n\in \mathbb{Z}$ such that $F_n\neq 0$. The full subcategory of $\dcat A$ consisting of perfect complexes is denoted by $\perf A$.

\subsection*{Derived actions}

Let $B$ be an $A$-algebra. We say that an $A$-complex $F$ admits a \emph{derived action} of $B$, or that $F$ is a \emph{derived $B$-complex}, to mean that there is a map of $A$-algebras
\[
B\longrightarrow \End_{\dcat A}(F)\,;
\]
equivalently, the canonical map $A\to \End_{\dcat A}(F)$ of rings factors through $B$. In this case, each $\hh_i(F)$ has the structure of a $B$-module, compatible with the $A$-action. However the existence of a derived $B$-action is a strictly stronger condition.

We write $\dcat[B]A$ for the subcategory (not necessarily full) of $\dcat A$ consisting of derived $B$-complexes; a morphism $f\colon F\to G$ in this category is a morphism in $\dcat A$ that is compatible with the $B$-actions, in the obvious way. In the same vein, $\perf[B]A$ are the derived $B$-complexes in $\perf A$. In what follows, the structure of these categories does not come into play; we only introduce them for the convenience of notation. In this language, Question~\ref{qu:Shekhar} has the more succinct formulation:

\begin{question}
   \label{qu:Shekhar2}
   Let $A$ be a noetherian local ring and $B$ be a noetherian local $A$-algebra. For any $A$-complex $F$ in $\perf[B]A$ with $\inf\hh_*(F)= 0$ and satisfying
   \[
   \pdim_AF\le \edim A - \edim B\,,
   \]
   the $B$-module $\hh_0(F)$ is free.
\end{question}

As noted above, a derived $B$-action on $F$ implies that the $A$-module structure on $\hh_*(F)$ extends to that of a $B$-module. When the $B$-module $\hh_0(F)$ is free, the $B$-module $\hh_*(F)$ is faithful and $\dim_A F = \dim B$. Note that even these statements need not to hold if the hypothesis in the question is weakened to asking only for the existence of a $B$-action on $\hh_*(F)$. Here is an example that illustrates this point; see however Theorem~\ref{th:regular_case}.

\begin{example}
\label{ex:H-action}
Let $A$ be a noetherian local ring for which there is an element $x\in \fm_A\setminus \fm_A^2$ that is in the socle of $A$ and such that $\fm_A\ne (x)$. Set $B=A/xA$. Let $f$ be an element in $\fm_A \setminus (x)$ and consider the complex
\[
F\colonequals 0\longrightarrow A^2 \xrightarrow{[x, f]} A\longrightarrow 0
\]
Then $\hh_0(F)=A/(x,f)$ and $x\cdot \hh_1(F)=0$, since $\hh_1(F)\subseteq \fm_AA^2$, as can be checked directly. Thus $\hh_*(F)$ is a $B$-module, and $\edim A - \edim B=1$. However $\hh_0(F)$ is not free as a $B$-module. For instance, take $A\colonequals k[x,y]/(x^2,xy)$ and $f=y$.
\end{example}

\subsection*{Exceptional complete intersection maps}
A local homomorphism $\vf\colon A\to B$ of noetherian local rings is a \emph{complete intersection} if, possibly after completing $B$ at its maximal ideal, $\vf$ factors as $A\to A'\to B$ where $A\to A'$ is flat map of noetherian local rings with regular closed fiber, and $A'\to B$ is surjective with kernel generated by a regular sequence. If in addition the regular sequence is linearly independent in $\fm_{A'}/\fm_{A'}^2$, then $\vf$ is an \emph{exceptional complete intersection}; see \cite{Brochard/Iyengar/Khare:2023a} for details.

\section{The regular case}
\label{section:regular_case}
In this section we prove the result below.

\begin{theorem}
 \label{th:regular_case}
  Let $A$ be a noetherian local ring, $B$ a noetherian local $A$-algebra, and $F$ a perfect $A$-complex with $\inf\hh_*(F)= 0$. Suppose that the $A$-module structure on $\hh_*(F)$ extends to that of a $B$-module, and that
  \[
  \pdim_AF \leq \edim A-\edim B\le	 \depth A - \dim B\,.
  \]
 The following conclusions hold:
 \begin{enumerate}[\quad\rm(1)]
     \item  The $B$-module $\hh_0(F)$ is free and $\sup\hh_*(F)=0$.
     \item  Both inequalities above are equalities.
     \item  The rings $A$ and $B$ are Cohen-Macaulay.
    \item The map $A\to B$ is an exceptional complete intersection.
 \end{enumerate} 
\end{theorem}

In the statement, the inequality displayed on the right is equivalent to
\[
\edim A - \depth A \le \edim B-\dim B\,,
\]
which roughly means that $A$ is ``more regular than $B$''. It holds if $A$ is regular, for the latter property is characterized by the condition that $\edim A = \depth A$. Hence the result generalizes \cite[Section 6]{Calegari/Geraghty:2018}. 

The proof of Theorem~\ref{th:regular_case} is given towards the end of this section. We extract a couple of computations that go into the argument.

\begin{lemma}
\label{le:NIT}
Let $A$ be a noetherian local ring, $B$ a noetherian local $A$-algebra, and $F$ a perfect $A$-complex  with $\inf\hh_*(F)= 0$. If the $A$-module structure on $\hh_*(F)$ extends to that of a $B$-module, then the following inequalities hold.
 \begin{enumerate}[{\quad\rm(1)}]
  \item   $\pdim_AF\geq \depth A-\dim B+\sup\hh_*(F)$;
  \item    $\pdim_AF \geq \dim A - \dim B$.
  \end{enumerate}
  Equality holds in \emph{(2)} if and only if the Cohen-Macaulay defect of $F$ equals that of $A$ and $\dim_A \hh_{0}(F)=\dim B$.
\end{lemma}

\begin{proof} 
Observe that if a $B$-module $M$ is finitely generated over $A$, then $\dim_B M = \dim_A M$. Indeed, the injective ring map $A/\Ann_A(M) \to B/\Ann_B(M)$ is finite, because $B/\Ann_B(M)$ is an $A$-submodule of $\End_A(M)$, which is a finite A-module. This applies in particular to all the modules $\hh_i(F)$.

Since $F$ is a perfect $A$-complex one has
\begin{equation}
\label{eq:abit}
\depth A - \depth_AF = \pdim_AF \ge \dim A - \dim_AF\,.    
\end{equation}
The equality is the Auslander-Buchsbaum formula; see~\cite[Proposition~1.3]{Iversen:1977}. The inequality is one form of the New Intersection Theorem proved by Roberts~\cite{Roberts:1987a}; see Iversen~\cite[Theorem~4.1]{Iversen:1977}. These estimates are used below.

(1) We can assume that $F$ is finite free and minimal. Set $p \colonequals \pdim_AF$ and $s\colonequals \sup \hh_*(F)$, and consider the complex
\[
 F' \colonequals 0\longrightarrow  F_p \longrightarrow  \dots \longrightarrow  F_s\longrightarrow  0 \,.
\]
This is a minimal free resolution of $\hh_s(F')$, of length $p -s$, so one gets the first of the following sequence of (in)equalities:
\begin{align*}
 p  &= s + \pdim_A\hh_s(F') \\
    &= s +\depth A  - \depth_A \hh_s(F') \\
    &\geq s +\depth A - \dim_A\hh_s(F) \\
    &= s + \depth A - \dim_B\hh_s(F) \\
    &\geq s + \depth A - \dim B
\end{align*}
The second equality is  by the Auslander-Buchsbaum formula, which applies because $\hh_s(F')$ has finite projective dimension; see \eqref{eq:abit}. The first inequality holds because $\hh_s(F)$ is a nonzero submodule of $\hh_s(F')$.  This justifies part (1).

(2) We have the following sequence of (in)equalities.
\begin{align*}
\dim_A F & = \sup\{\dim_A \hh_n(F) - n\mid n\in \mathbb{Z}\} \\    
         & = \sup\{\dim_B \hh_n(F) - n\mid n\in \mathbb{Z}\} \\    
         &\leq \sup\{\dim B - n\mid n\in \mathbb{Z} \text{ and } \hh_n(F)\ne 0\} \\    
         &= \dim B\,.
\end{align*}
The first equality is the definition of dimension and the inequality is clear. Given these and the inequality in \eqref{eq:abit}, part (2) follows, as does the claim about when equality holds.
\end{proof}

\begin{proof}[Proof of Theorem~\ref{th:regular_case}]
The hypotheses and Lemma~\ref{le:NIT}(1) imply (2) and also that $\sup\hh_*(F)=0$, that is to say, $F$ is a free resolution of $\hh_0(F)$. Hence 
\[
\pdim_A \hh_0(F)=\edim A-\edim B\,,
\]
so~\cite[Theorem 3.1]{Brochard/Iyengar/Khare:2023a} yields that the $B$-module $\hh_0(F)$ is free, and also (4).

Moreover, Lemma~\ref{le:NIT}(2) and the hypothesis yields $\dim A\le \depth A$, so the ring~$A$ is Cohen-Macaulay, and since the map $A\to B$ is complete intersection, it follows that $B$ is Cohen-Macaulay as well; see, for instance, \cite[Section~2.1]{Bruns/Herzog:1998a}.
\end{proof}

\section{The case of a complex of \texorpdfstring{$B$}{B}-modules}
\label{section:B_modules}

In this section we prove the result below, which is Theorem~\ref{th:main1}.  

\begin{theorem}
\label{th:B-module_case}
 Let $\vf\colon A\to B$ be a local map of noetherian local rings and $F$ a complex of $B$-modules with $\hh_*(F)$ finitely generated over $B$,  $\inf \hh_*(F)=0$, and 
 \[
 \fdim_AF\le \edim A-\edim B\,.
 \]
 Then the $B$-module $\hh_0(F)$ is free, $\hh_i(F)=0$ for $i\ne 0$, and the map $\vf$ is an exceptional complete intersection. 
\end{theorem}

For the notion of flat dimension for complexes see, for instance, \cite{Avramov/Foxby:1991a}. Given a local homomorphism $\vf\colon A\to B$ of local rings and a $B$-complex $F$ such that the $B$-module $\hh_*(F)$ is finitely generated, the flat dimension of $F$ over $A$ is finite if and only if $\Tor^A_i(k_A,F)=0$ for $i\gg 0$; see \cite[Proposition~5.5]{Avramov/Foxby:1991a}.

As explained in the Introduction, one can view this statement as providing an obstruction to the existence of a lift of the derived $B$-action on $F$ to a strict $B$-action; said otherwise, to lifting the given map of $A$-algebras $B\to \End_{\dcat{A}}(F)$ to a map of dg  (= differential graded) $A$-algebras $B\to \End_A(F)$, where $\End_A(F)$ is the endomorphism dg $A$-algebra of $F$.

 The proof of Theorem~\ref{th:B-module_case} uses dg algebra methods. In what follows, by a dg algebra we mean one that is strictly graded-commutative. For basic definitions and constructions regarding dg algebras and dg modules, see \cite{Avramov:1999a}. In fact, the only dg algebra we have to work with is the Koszul complex on a finite set of elements.

\begin{example}
\label{ex:koszul_as_a_dg_algebra}
 Let $\bs a\colonequals a_1,\dots, a_c$ be elements in a commutative ring $A$. The \emph{Koszul complex} on $\bs a$, denoted $K(\bs a)$, is the exterior algebra on a free module $\oplus_{i=1}^c Ae_i$, with $|e_i|=1$. Endow it with a differential that maps $e_i$ to $a_i$, and satisfies the Leibniz rule. Then $K$ is a dg algebra over $A$.
\end{example}

For homological invariants, like Tor, of dg modules see \cite{Avramov:1998a, Avramov/Buchweitz/Iyengar/Miller:2010a}. By a local dg algebra $R$ we mean that $R_0$ is a local ring and the $\hh_0(R)$-module $\hh_i(R)$ is finitely generated for each $i$. We write $k_R$ for the residue field of $R_0$. The Koszul dg algebra on a finite set of elements over a local ring is the basic example of a local dg algebra.

\begin{chunk}
\label{ch:Betti}
 Let $R$ be a local dg algebra and $M$ a dg $R$-module such that the $\hh_0(R)$-module $\hh_*(M)$ is finitely generated.  Extending the notions from modules over local rings, the \emph{Betti numbers} of $M$ are the numbers
 \[
 \beta_n^R(M)=\dim_{k_R}\Tor^R_n(k_R,M)\,,
 \]
and their generating series is the \emph{Poincar\'e series} of $M$, namely
\[
 P_M^R(t)=\sum_{n\in \mbb{Z}} \beta_n^R(M)t^n.
\]
 The \emph{projective dimension} of $M$ over $R$ is
 \[
  \pdim_RM\colonequals \sup\{n\in \mbb{Z}\mid \beta^R_n(M)\ne 0\}\,.
 \]
 \end{chunk}

Here is a criterion for detecting free dg modules.

\begin{lemma}
\label{le:free}
Let $R$ be a local dg algebra and let $M$ be a dg $R$-module such that the $\hh_0(R)$-module $\hh_*(M)$ is finitely generated and $\inf\hh_*(M)=0$. Then $M\simeq R^b$ in $\dcat R$, for some integer $b$, if and only if $P_M^R(t)=b$.
\end{lemma}

\begin{proof}
The ``only if" condition is obvious. Conversely, assume that $P_M^R(t)=b$, i.e. that $\beta^R_0(M)=b$ and $\beta_n^R(M)=0$ for $n\geq 1$. Writing $k$ for the residue field of $R_0$, the right exactness of tensor products yields $\Tor^R_0(k,M)\cong k\otimes_{R_0}\hh_0(M)$, so that $b=\nu_{R_0}\hh_0(M)$. Pick cycles $z_1,\dots,z_b$ in $M_0$ whose residue classes in $\hh_0(M)$ form a minimal generating set for it, as an $R_0$-module. Consider the dg $R$-module $\oplus_{i=1}^b Re_i$ where $d(e_i)=0$ for each $i$ and consider the morphism $f\colon \oplus_{i=1}^b Re_i\to M$ where $e_i$ maps to $z_i$. Complete it to an exact triangle 
\[
\bigoplus_{i=1}^b Re_i \xrightarrow{\ f\ }  M \longrightarrow W \longrightarrow
\]
in $\dcat A$. It is immediate from the corresponding long exact sequence in homology that the $\hh_0(R)$-module $\hh_*(W)$ is finitely generated.  Moreover, the construction of $f$ and the hypothesis on $M$ yield that $\Tor^R_*(k,W)=0$. Hence $W=0$: else for  $i=\inf\hh_*(W)$, one has $\Tor^R_i(k,W)\cong k\otimes_{R_0}\hh_i(W)\ne  0$, again by the right exactness of tensor products and Nakayama's Lemma. In conclusion, $f$ is a quasi-isomorphism, as desired.
\end{proof}

The following proposition, which is a generalization of Nagata's theorem on surjective exceptional complete intersection maps, see~\cite[Proposition~3.3.5(1)]{Avramov:1998a}, is the key input in the proof of Theorem~\ref{th:B-module_case}.

\begin{proposition}
\label{pr:large_map}
 Let $A$ be a noetherian local ring, and fix elements $\bs a\colonequals a_1,\dots, a_c$ in $\fm_A$ such that the image of the set $\bs a$ in $\fm_A/\fm_A^2$ is linearly independent. Let $K$ be the Koszul dg $A$-algebra on $\bs a$. For any dg $K$-module $M$ such that the $\hh_0(K)$-module $\hh_*(M)$ is finitely generated there is an equality
\[
 P_M^A(t)=P_M^K(t).(1+t)^c.
\]
In particular, $\pdim_AM$ and $\pdim_KM$ are finite simultaneously, and
\[
\pdim_AM=\pdim_KM+c\,.
\]
\end{proposition}

\begin{proof}
Write $k$ for the residue field of $A$ and consider the graded $k$-vector space
\[
V\colonequals \bigoplus_{n\geqslant 0} \Sigma^n k^{\binom cn}\,.
\]
viewed as a dg $K$-module via the augmentation $K\to k$. In the derived category of dg $K$-modules, one has a quasi-isomorphism $k\otimes_AK \simeq V$. This can be verified by arguing exactly as in the proof of \cite[Theorem~9.1]{Avramov/Buchweitz/Iyengar/Miller:2010a}. This yields quasi-isomorphisms
\[
k\lotimes_A M 
    \simeq (k\lotimes_A K)\lotimes_K M
     \simeq V \lotimes_k (k\lotimes_KM)\,.    
\]
In homology this yields an isomorphism
\[
\Tor^A_*(k,M)\cong \hh_*(V)\otimes_k \Tor^K_*(k,M)\,.
\]
This yields the stated equality of Poincar\'e series. The rest of the assertions are immediate consequences of this equality.
\end{proof}

\begin{proof}[Proof of Theorem~\ref{th:B-module_case}]
First we reduce to the case where the local ring $B$ is complete at its maximal ideal $\fm_B$, as follows.

Let $\iota\colon B\to \hat B$ denote the completion of $B$ at  $\fm_B$, and set $\hat F= \hat B \otimes_B F$. It is immediate from the definition that $\vf$ is an exceptional complete intersection if and only if $\iota\vf$ is an exceptional complete intersection. Moreover, since the map $B\to \hat B$ is faithfully flat, one has
 \[
 \fdim_A \hat F =\fdim_A F\quad \text{and}\quad \hh_i(\hat F)\cong \hat B\otimes_B \hh_i(F)
 \]
 for all $i$. Hence the $\hat B$-module $\hh_*(\hat F)$ is finitely generated, and $\hh_i(F)=0$ for $i\ne 0$ if and only if $\hh_i(\hat F)\ne 0$. Moreover, the $B$-module $\hh_0(F)$ is free if and only if the $\hat B$-module $\hh_0(\hat F)$ is free. Thus, replacing $\vf$ by $\iota\vf$ and $F$ by $\hat F$ we can assume $B$ is $\fm_B$-adically complete.

Next we reduce to the case where $\vf$ is surjective. Since $B$ is complete, $\vf$ admits a factorization  $A\to A'\to B$ where $A\to A'$ is a flat local map whose closed fiber $A'/\fm_A A'$ is regular, and the map $A'\to B$ is surjective; see  \cite[1.1]{Avramov/Foxby/Herzog:1994}.
We claim that
\begin{align*}
    \fdim_{A'}F 
    & \leq \fdim_AF+\edim(A'/\fm_AA') \\
    &=\fdim_AF+\edim A' -\edim A\\
    &\le \edim A - \edim B + \edim A' - \edim A\\
    &=\edim A'-\edim B\,.
\end{align*}
The first inequality can be verified by mimicking the proof of the inequality given in \cite[2.7]{Avramov/Foxby/Halperin:descent_and_ascent}  dealing with the case where $F$ is a module. For the equality see~\cite[2.6]{Brochard/Iyengar/Khare:2023a}, whereas the second inequality is by our hypotheses. Finally $\vf$ is complete intersection if and only if the map $A'\to B$ is complete intersection, again by definition. Thus  replacing $\vf$ by the map $A'\to B$ brings us to the context of surjective maps.

Since $\vf$ is surjective  $\hh_*(F)$ is finitely generated also as an $A$-module, and hence $\fdim_{A}F=\pdim_{A}F$.  Let $c\colonequals \edim A -\edim B$. Since $\varphi$ is surjective, there exists a sequence $\bs a\colonequals a_1,\dots, a_c$ of elements in $\ker(\varphi)$ such that the image of the set $\bs a$ in $\fm_A/\fm_A^2$ is linearly independent. Let $K$ be the Koszul dg $A$-algebra  on $\bs a$. Since $\hh_0(K)=A/(\bs a)$ one gets morphisms of dg algebras $ K\to \hh_0(K)\to B$. Via the morphism $K\to B$ we view $F$ as a dg module over $K$. Since $F$ has finite projective dimension over $A$, Proposition~\ref{pr:large_map} yields that it has finite projective also over $K$, and
\begin{align*}
\pdim_KF & = \pdim_A F -c \\
         & \le \edim A-\edim B - c\\
         & \le 0\,.
\end{align*}
Thus $\beta^K_i(F)=0$ for $i\ge 1$. Since $\inf\hh_*(F)=0$ and $\hh_*(F)\ne 0$ we deduce that $F\simeq  K^b$ as dg $K$-modules for some integer $b\ge 1$; see Lemma~\ref{le:free}.  In particular $\hh_0(F)$ is a direct sum of copies of $\hh_0(K)$. Since it is also a $B$-module this implies that the kernel of $A\to B$ is $(\bs a)$, hence $B=\hh_0(K)$ and $\hh_0(F)$ is free as a $B$-module.
The last conclusion implies that in $\dcat B$ there is a quasi-isomorphism
\[
F\simeq \tau_{\geqslant 1}(F)\oplus \hh_0(F)\,.
\]
Since the action of $K$ on $F$ factors through $B$, this is also a quasi-isomorphism of dg $K$-modules. In particular, $\pdim_K F$ finite implies $\pdim_K\hh_0(F)$, and hence also $\pdim_KB$, is finite. Now the amplitude inequality~\cite[Theorem~0.2]{Jorgensen:2010} applied to the dg $K$-algebra $K$ yields that
$\inf\hh_*(K)\le 0$, that is to say, $\hh_i(K)=0$ for $i\ne 0$; equivalently, the natural map $K\to B$ is a quasi-isomorphism. This implies that $\bs a$ is a regular sequence in $A$. This proves everything we need.
\end{proof}

\begin{example}
\label{ex:koszul}
Let $A$ be a noetherian local ring and $\bs x\colonequals x_1,\dots,x_p$ elements that form part of a minimal generating set for the ideal $\fm_A$. Set $B\colonequals A/(\bs x)$, so $\edim A - \edim B = p$.  The natural map $A\to \End_{\cat D}(K(\bs x))$ factors through $B$; see \cite[Proposition~1.6.5]{Bruns/Herzog:1998a}. Thus, the hypothesis of Question~\ref{qu:Shekhar} holds for the $A$-complex $K(\bs x)$, as does the conclusion, for $\hh_0(K) = B$. 

However, one can choose $A$ such that $\fm_A$ contains only zero-divisors, also the sequence $\bs x$ such that $\hh_i(K)$ is not free for $i\ge 1$. For example, $A\colonequals k[x,y]/(x^2,xy)$ and the element $x$. Then $\hh_1(K) \cong \fm_A$.
\end{example}

\section{A consequence for patching}
\label{section:nt}
In this section we spell out a consequence of Theorem~\ref{th:B-module_case} in number theory.  Consider the complexes $C_K$ that occur in \cite[\S 13]{Iyengar/Khare/Manning:2024b}.   We freely use the notation and set-up  there. So  let $F$ be a number field and consider the algebraic group $\PGL_{2}$ over $F$. Let $K_\infty\subseteq \PGL_{2}(F\otimes_{\mbb{Q}}{\mbb{R}})$ be a maximal compact subgroup. For any compact open subgroup  $K = \prod_vK_v \subseteq \PGL_2(\mbb{A}_{F}^\infty)$  consider the topological space:
\[
Y_K = \PGL_2(F)\backslash \PGL_2(\mbb{A}_F)/KK_\infty\,.
\]
It is an orbifold of dimension $2r_1+3r_2$. We say that $K$ is \emph{sufficiently small} if $gKg^{-1} \cap \PGL_2(F)$ is  torsion-free for all $g \in \PGL_2(\mbb{A}_F^\infty)$. We consider only such $K$. 

For any  such $K$, let $C_{K}$ be the complex of singular chains on $Y_K$ with coefficients in $\mcO$,   the ring of integers of a finite extension of $\mbb{Q}_p$, so that $C_{K}$ computes the homology $\hh_*(Y_K,\mcO)$. Then $C_{K}$ is quasi-isomorphic to a bounded complex of free $\mcO$-modules, and  so may be viewed as a perfect complex in $\dcat{\mcO}$.   We have a derived action of the Hecke algebra on $C_K$, and thus we may consider it as a subalgebra of  $ \End_{\dcat{\mcO}}(C_{K})$.   It is natural to ask if this Hecke action can be lifted to  $ \End_{{\mcO}}(C_{K})$;    Theorem~\ref{th:B-module_case} answers  a version of this question that arises in patching arguments   in the negative.

Namely  consider  a level subgroup $K=K_{\Sigma, Q}$ defined in  loc. cit.,  and a maximal ideal $\fm$  of the Hecke algebra  in $ \End_{\dcat{\mcO}}(C_{K})$,  with corresponding Galois representation $\rhobar_{\fm}$ that  satisfies the Taylor-Wiles hypothesis,  and with  $Q$  a set of Taylor-Wiles primes for $\rhobar_{\fm}$  as in \cite[Proposition 15.2]{Iyengar/Khare/Manning:2024b}. Consider the complex  $C_{\Sigma, Q}$ in   \cite[\S  13]{Iyengar/Khare/Manning:2024b}. Thus $C_{\Sigma, Q}$  is a perfect complex of $A=\mcO[\Delta_Q]$-modules, with $Q$ a set of Taylor-Wiles primes in Proposition 15.2 of loc. cit.,  such that the map $A \to  \End_{\dcat{A}}(C_{\Sigma, Q})$ factors through a Hecke algebra $\mbb{T}_{\Sigma,Q}$ defined there.   Here $\Delta_Q$ is the  maximal $p$-quotient  of ${(\mcO_F/\prod_{\mathfrak{q}\in Q} \mathfrak{q})}^*$.

\begin{corollary}
With notation as above, assume Conjectures (A), (B), (C) of \cite{Iyengar/Khare/Manning:2024b} and that the field $F$ is not totally real. If
$\hh_*(Y_{K_{\Sigma,Q}},\mcO)_{\fm} \otimes_\mcO \mbb{Q}$ is non-trivial, then  $C_{\Sigma, Q}$ is not quasi-isomorphic  in $\dcat{\mcO[\Delta_Q]}$ to a complex of $\mbb{T}_{\Sigma,Q}$-modules.
\end{corollary}

We can also deduce  that there is no perfect complex of  $\mcO[\Delta_Q]$-modules   with  a Hecke action and which computes the Hecke module $ \hh_*(Y_{K_{\Sigma,Q}},\mcO)$.

\begin{proof}
 As explained above, for $A=\mcO[\Delta_Q]$ and $B=\mbb{T}_{\Sigma,Q}$, the $A$-complex $C_{\Sigma,Q}$ is finite free and admits a derived $B$-action.  Conjectures (A) and (B) imply that there is  surjective map  $R_{\Sigma,Q} \to B$. Conjecture (C) and the choice of Taylor-Wiles primes $Q$ coming from Proposition 15.2 of \cite{Iyengar/Khare/Manning:2024b}   imply that there is an    inequality 
\[
\pdim_AC_{\Sigma,Q}  \le  \edim A - \edim R_{\Sigma,Q} \,,
\]  
and hence one gets an inequality
\[
\pdim_AC_{\Sigma,Q}  \le  \edim A - \edim B \,.
\]
We are thus in the setting of Theorem~\ref{th:B-module_case}, where the $A$ and $B$ are as above, and the complex in question is $C_{\Sigma,Q}$. Because of our assumption that    $\hh_*(C_{\Sigma,Q}) \otimes {\mbb{Q}}$ is non-trivial, and that the field $F$ is not totally real,  we also have that $\hh_*(C_{\Sigma,Q}) \otimes {\mbb{Q}}$ is  not concentrated in a single degree; see, for instance, \cite[Theorem III.5.1]{Borel/Wallach:2000}.  Thus Theorem~\ref{th:B-module_case} gives the desired conclusion.
\end{proof}

In the corollary, when  $F$ is an imaginary quadratic field, we need to only  assume  Conjectures (A) and (B)  which are certain `local-global compatibility’ statements for Galois representations $\rho\colon G_F \to \mathrm{GL}_2(\mbb{T}_{\Sigma,Q})$ (and that are known up to going modulo a  nilpotent ideal).

\subsection*{Poincar\'e duality}
Let $A,B$ and $F$ be as in Question~\ref{qu:Shekhar}.  In the number theoretic context (as in this section) the complex $F$ has Poincar\'e duality compatible with the derived $B$-action. By this we mean that there is an isomorphism $\RHom_A(F,\omega_A)\simeq \Sigma^n{F}$ in $\dcat[B]{A}$. Here $\omega_A$ is the dualizing complex of $A$, normalized as in \cite{StacksProject}; see also \cite[Section~4]{Iyengar/Khare/Manning:2024b}. So it is natural to wonder if this has a bearing on Question~\ref{qu:Shekhar}. 

There is some evidence pointing in this direction: An affirmative answer to Question~\ref{qu:Shekhar} implies that the $B$-module $\hh_*(F)$ is faithful. In  \cite[Lemma~10.8]{Iyengar/Khare/Manning:2024b} it is proved that when $p\le 1$ and $F$ has Poincar\'e duality, $\hh_*(F)$ is faithful as a module over all of $\End_{\dcat A}(F)$;  this does not even require $F$ to have a derived $B$-action.

\section{The complete intersection case}
\label{section:ci_case}
In this section we prove Theorem~\ref{th:main2}; see Corollary~\ref{co:ci-fiber} below. 

Let  $(A,\fm_A,k_A)$ be a noetherian local ring and $M$ an $A$-complex such that the $A$-module $\hh_*(A)$ is finitely generated. The Betti numbers of $M$ are as in \ref{ch:Betti}.  For $a=\inf\hh_*(M)$ the right-exactness of the tensor product yields an isomorphism of $k_A$-vector spaces
 \[
 \Tor^A_a(k_A,M)\cong k_A\otimes_A\hh_a(M)\,.
 \]
Hence $\beta^A_a(M)$ is the minimal number of generators of  $\hh_a(M)$. 
 
In what follows, given a sequence of elements $\bs x$ in $A$, we write $K(\bs x)$ for the Koszul complex on $\bs x$ with coefficients in $A$; see \cite[\S1.6]{Bruns/Herzog:1998a}, and also Example~\ref{ex:koszul_as_a_dg_algebra}.

\begin{theorem}
\label{th:ci_case}
 Let $A$ be a noetherian local ring, $B$ a finite local $A$-algebra, and $F$ in $\perf[B]A$ with $\inf \hh_*(F) = 0$ such that
\[
p\colonequals \pdim_A F \leq \edim A -  \beta^B_0(\fm_AB)\,.
\]
The following statements hold:
\begin{enumerate}[\quad\rm(1)]
 \item The $B$-module  $M=\hh_0(F)$ is free and $\beta^A_1(M) = p \cdot \beta^A_0(M)$.

\item The ring $B/\fm_AB$ is a complete intersection of dimension  zero.

\item One has $p = \edim A - \beta^B_0(\fm_AB) = \edim A-\edim B$. 

\item 
The Betti numbers of $F$ satisfy inequalities
\[
 \beta^A_i(F) =\binom{p}{i} \beta^A_0(F) \quad\text{for} \quad 0\leq i\leq p\,.
\]
\end{enumerate}
 If moreover the kernel of the map $A\to B$ contains a set $\bs x$ of $p$ elements that extend to a minimal generating set for the ideal $\fm_A$, then $F$ is quasi-isomorphic to a direct sum of copies of the Koszul complex $K(\bs x)$.
\end{theorem}

The following corollary justifies Theorem~\ref{th:main2}.

\begin{corollary}
\label{co:ci-fiber}
Question~\ref{qu:Shekhar} has a positive answer when the ring $B/\fm_AB$ is a complete intersection.    
\end{corollary}

\begin{proof}
 Since $B/\fm_A B$ is a complete intersection one has $\beta^B_0(\fm_AB)\le \edim B$, by~\cite[5.4]{Brochard/Mezard:2011}. Thus Theorem~\ref{th:ci_case} applies and yields the desired conclusion.
\end{proof}

\begin{proof}[Proof of Theorem~\ref{th:ci_case}]
Replacing $F$ by its minimal free resolution over $A$, we can assume that $F$ is a finite free $A$-complex satisfying $d(F)\subseteq \fm_AF$. It suffices to verify statements (3) and (4) because the latter gives the second inequality below
\[
\beta^A_1(M) \leq \beta^A_1(F) \leq p\cdot \beta^A_0(F) =p\cdot \beta^A_0(M)\,.
\]
The inequality on the left and the equality hold because $F_1\to F_0\to M$ is a presentation of $M$ with $d(F_1)\subseteq \fm_A F_0$. Thus (1) and (2) follow by~\cite[4.12]{Brochard:2023}.

(3) and (4): Since $\beta^B_0(\fm_AB)\leq \edim A-p$  there exist a minimal system $x_1,\dots, x_n$ of generators of $\fm_A$ such  that  $x_{p+1}, \dots, x_n$ generate $\fm_AB$.  Since $F$ is  minimal the differential $d\colon F \to F$ can be written $d =\sum_{j=1}^nx_jd_{j}$.

With $\varphi\colon A\to B$ the structure map, for any $i\leq p$ there is a  relation
\[
\varphi(x_i)=\sum_{j=p +1}^{n}\varphi(x_j)b_{ij} \quad\text{with $b_{ij}\in B$.}
\]
We write $\cat D$ for the derived category and $\cat K$ for the homotopy category of $A$, respectively.  Since $F$ is a complex of free $A$-modules, $ \End_{\cat D}(F)\cong \End_{\cat K}(F)$. Noting that $\End_{\cat K}(F)$ is the homology of the complex $\End_A(F)$, for any  $i, j$ we choose a lift $\widetilde{b_{ij}}$ in  $\End_{A}(F)$ of the image of $b_{ij}$ in  $\End_{\cat K}(F)$. Then  $x_i\id_{F}-\sum_{j=p +1}^{n} x_j\widetilde{b_{ij}}$ is homotopic to $0$; let $h_i$ be a homotopy witnessing this, so that
\begin{equation}
\label{eq:homotopy}
 x_i\id_{F}-\sum_{j=p +1}^{n}x_j\widetilde{b_{ij}} = d h_{i}+h_{i}d = \sum_{j=1}^n x_j(d_{j}h_{i}+h_{i}d_{j})
\end{equation}
$1\leq i\leq p$. Since the complex $\End_A(F)$ is $A$-free and $x_1,\dots,x_n$ is a minimal generating set for the ideal $\fm_A$, it follows that for any $i,j$ with $1\leq  i,j\leq p $ we have
\begin{alignat*}{3}
\overline{d}_{j}\overline{h}_{i}+\overline{h}_{i}\overline{d}_{j}&=0
 &&\textrm{if } i\neq j, \textrm{and} \\
\overline{d}_{j}\overline{h}_{i}+\overline{h}_{i}\overline{d}_{j}&=\id_{\overline{F}} \quad  &&\textrm{if } i=j,
\end{alignat*}
in $\End_{k_A}(\overline{F})$, where the bar denotes reduction modulo $\fm_A$. 

Let $W$ be the free associative $k_A$-algebra on variables $s_1,\dots,s_p,t_1,\dots,t_p$ modulo the ideal generated by indeterminates $\{s_it_j+t_js_i=\delta_{ij}\}_{1\le i,j\le p}$, where $\delta_{ii}=1$ and $\delta_{ij}=0$ for $i\ne j$. We view this as a graded $k_A$-algebra with $|s_i|=1$ and $|t_i|=-1$ for all $i$. It follows from the relations above that $\overline{F}$ has the structure of a graded left $W$-module, with $s_i$ acting via $\overline{h}_i$ and $t_i$ acting via $\overline{d}_i$. Since $\overline{F}_0\ne 0$ and $\overline{F}_n=0$ for $n\ge p+1$, one can apply Proposition~\ref{pr:Wmod} to deduce that the natural map $W\otimes_{k_A}\overline{F}_0\to F$ induces an isomorphism
\[
\phi\colon k_A[\bs s] \otimes_{k}\overline{F_0}\xrightarrow{\ \cong\ } \overline{F}
\]
of graded left $W$-modules, where $k_A[\bs s]$ is the exterior algebra on $\bs s$, viewed as a quotient of $W$; see Appendix~\ref{se:appendix} and Definition~\ref{de:exterior-alg}. This proves~(3) and (4). 

The map $\phi$ can be described explicitly and this is useful in what follows: For any subsequence $I= (i_1< \dots < i_n)$ of $[p]=(1<\cdots <p)$, set
\[
s_I= s_{i_1}\cdots s_{i_n}\quad \text{and}\quad {h}_I = {h}_{i_1}\cdots {h}_{i_n}\,.
\]
Then $k_A[\bs s]$ has basis $\{s_I\}_{I\subseteq [p]}$, and by construction, the map $\phi$ above is given by
\[
s_I\otimes f\mapsto  \overline{h}_I(f) \quad \text{for $f\in \overline{F_0}$.}
\]

Suppose $\ker(\varphi)$ contains $p$ elements that extend to a minimal generating set for the ideal $\fm_A$. Then, in the notation above, we may assume  $\bs x\colonequals x_1,\dots, x_p$ map to 0 in $B$.  We can choose $\widetilde{b_{ij}}=0$ so that the  relation~\eqref{eq:homotopy} becomes
\begin{equation}
\label{eq:homotopy2}
  x_i\id_{F}=d h_i+h_id\,.
\end{equation}
Let $K$ be the Koszul complex on $\bs x$. Then $K_1=\oplus_{i=1}^p Ae_i$ with $d(e_i)=x_i$, and as a graded free $A$-module
\[
K = \bigoplus_{I\subseteq [p]} Ae_I\,,
\]
where $e_I$ has the obvious meaning. Consider the map of graded $A$-modules
\[
\Phi \colon K\otimes_A F_0\longrightarrow F \quad\text{where $e_I\otimes f\mapsto h_I(f)$.}
\]
It is straightforward to verify, by an induction on $|I|$ and using \eqref{eq:homotopy2}, that $\Phi$ is also a morphism of complexes.  Finally, observe that $k_A\otimes_A\Phi =\phi$. Since $\phi$ is an isomorphism so is $\Phi$, being a map between free $A$-modules.
\end{proof}

Theorem~\ref{th:ci_case}, and more precisely, a step in its proof, yields a characterization of perfect complexes that are isomorphic to a direct sum of  Koszul complexes.

\begin{corollary}
\label{co:koszul-char}
 A perfect complex $F\colon 0\to F_p\to\cdots\to F_0\to 0$ over a noetherian local ring $A$ is quasi-isomorphic to a direct sum of the Koszul complex $K(\bs x)$, where  $\bs x\colonequals x_1,\dots,x_p$ is part of a minimal generating set for $\fm_A$, if and only if the kernel of the natural map $A\to \End_{\dcat A}(F)$ contains $\bs x$.
\end{corollary}

\begin{proof}
For any sequence $\bs x$ of length $p$, the Koszul complex $K(\bs x)$ is concentrated in degrees $0,\dots, p$, and the kernel of the map $A\to \End_{\dcat A}(K(\bs x))$ equals the ideal $(\bs x)$; see \cite[Proposition~1.6.5]{Bruns/Herzog:1998a}. This implies the only if direction of the corollary.

Suppose that $F$ is a perfect $A$-complex as above such that the kernel of the map $A\to \End_{\dcat A}(F)$ contains a subset $\bs x$ of length $p$ that extends to a minimal generating set for the ideal $\fm_A$. Set $B=A/(\bs x)$; then $F$ is a derived $B$-complex, and the hypothesis on $\bs x$ implies $p=\edim A - \edim B$. Since $\fm_A B=\fm_B$ Theorem~\ref{th:ci_case} applies and yields the desired conclusion.
\end{proof}

 We end the section with some examples of complexes that satisfy the assumptions of Question~\ref{qu:Shekhar}, which might help to understand the structure of the complex~$F$.

\begin{example}
\label{ex:not_koszul}
 Given Theorem~\ref{th:ci_case}(4) it might be tempting to believe that the last conclusion of \emph{op.~cit.} holds even without the additional assumption on the kernel of $\varphi$.  However, this is not the case: let $A=k[x,y]/(x,y)^2$ and $B=k[u]/(u^4)$. Consider the morphism $\varphi\colon A\to B$ defined by $\varphi(x)=u^2$ and $\varphi(y)=u^3$ and let $F$ be the complex 
 \[
 0\longrightarrow A^2 \xrightarrow{
 \begin{pmatrix}
    y & 0 \\ -x & y
   \end{pmatrix} }
   A^2\longrightarrow 0
 \]
 concentrated in degrees $0$ and $1$. Then $F$ satisfies the assumptions of the above theorem: to get a  $B$-action on $F$ we map $u$ to the morphism of complexes 
\[
U\colon F\to F \quad\text{defined by}\quad 
U_{0} = \begin{pmatrix}
             0 & x \\ 1 & 0
            \end{pmatrix} = U_1\,.
\]
Then $U^2=x\id_{F}$ and $U^3$ is homotopic to $y\id_{F}$; a homotopy is $-\id_{A^2}$.
Thus, the assignment $u\mapsto U$ indeed defines an $A$-algebra map from $B$ to $\End_{\dcat A}(F)$. 

However, $F$ is not isomorphic to a direct sum of Koszul complexes, even in the derived category of $A$, essentially because the matrix $d$ is not equivalent to a diagonal matrix.
\end{example}

\begin{example}
 Let us consider a similar example, but with $\edim A-\edim B=2$. Let $A=k[x,y,z]/(x,y,z)^2$ and $B=k[u]/(u^6)$ and $\varphi\colon A\to B$ the map of $k$-algebras that assigns $x, y$ and $z$ respectively to $u^3, u^4$ and $u^5$. Consider the following $A$-linear endomorphisms of $A^3$:
\begin{gather*}
    c_1\colonequals x E - y\id \quad\text{and} \quad    c_2\colonequals yE - z\id\,\quad\text{where}\\
    E =\begin{pmatrix}
       0&0&x \\ 1&0&0 \\ 0&1&0
      \end{pmatrix} \colon A^3\longrightarrow A^3\,.
\end{gather*}
Let $F$ be the $A$-complex
\[
0\longrightarrow A^3 \xrightarrow{\begin{pmatrix}  c_2 \\ - c_1 \end{pmatrix}} \begin{matrix} A^3 \\ \bigoplus \\ A^3 \end{matrix} 
                \xrightarrow{\begin{pmatrix} c_1 & c_2 \end{pmatrix}} A^3\longrightarrow 0
\]
This has a derived $B$-action where the element $u$ acts via the  endomorphism $U$ of $F$ prescribed by
\[
U_2 = E\,, U_1 = E\oplus E\,,\quad\text{and}\quad U_0=E\,.
\]
It is a straightforward verification to check that $U^3=x\id$, that $U^4$ is homotopic to $y\id$ and that $U^5$ is homotopic to $z\id$. So this defines indeed a derived action of $B$ on $F$. The cokernel $\hh_0(F)$ is a free $B$-module of rank 1.

 This example is essentially the only one. Indeed, let $F'$ be another minimal length~2 complex of free $A$-modules with a derived action of $B$. By Theorem~\ref{th:ci_case} we already know that $\hh_0(F')$ is a free $B$-module and that $\rank_A(F'_2)=\rank_A(F'_0)$ and $\rank_A(F'_1)=2\rank_A(F'_0)$. With a little more computation, we can prove that the complex is necessarily isomorphic to a direct sum of copies of the complex $F$ described above (the derived action of~$B$ is not unique however: there are other possibilities for $U_1$ and $U_2$, leading to non-isomorphic actions on $F'$ ; this does not change the $B$-module structure on~$\hh_0(F')$). 
 
 This proves in particular that $F'$ is isomorphic to a Koszul complex in the following sense. Let $A'=A[\![z_1, z_2]\!]$ and let $K$ be the Koszul complex on the sequence $z_1, z_2$ over $A'$. Let us provide $F'_0$ with a structure of $A'$-module, by making $z_1$ and $z_2$ act via $xu'_0-y\id$ and $yu'_0-z\id$ respectively. Then, as an $A$-complex, $F$ is isomorphic to $K\otimes_{A'}F'_0$.
\end{example}

 This example together with Theorem~\ref{th:ci_case}(4), and also Theorem~\ref{th:B-module_case}, suggest the following question.

\begin{question}
 In the context of Question~\ref{qu:Shekhar}, is there  a structure of $A'$-module on $F_0$, with $A'=A[\![z_1,\dots, z_c]\!]$ where $c=\edim(A)-\edim(B)$, such that $F\simeq K\otimes_{A'}F_0$, where $K$ is the Koszul complex on the sequence $z_1, \dots, z_c$ over~$A'$?
\end{question}

Under the assumptions of Theorem~\ref{th:ci_case}, when $F$ is minimal, by~(1) we must have $\ker(F_1\to F_0)\subseteq \fm_A F_1$. This is no longer true without the assumption that $B/\fm_AB$ is a complete intersection, as we see in the following example.

\begin{example}\rm
Let $A=k[x,y,z]/(x,y,z)^2$ and let $B=k[u,v]/(u,v)^4$. Let $\varphi\colon  A\to B$ be the morphism that maps $x, y$ and $z$ to $u^2, uv$ and $v^2$ respectively. Let $d\colon F_1\to F_0$ be a length complex of free $A$-modules, equipped with a derived action of $B$. Assume moreover that $H_0(F)$ is a free $B$-module (none of our theorems applies here, hence question~\ref{qu:Shekhar} is still open for this example). Then we can prove that $F$ is necessarily isomorphic to a direct sum of copies of the complex where $F_1=F_0=A^3$, and the differential $d$ is given by the matrix
\[
 [d]= \begin{pmatrix}
    -y & -z & 0 \\
    x&y&0 \\
    0&0&0
  \end{pmatrix}.
\]
The derived action of $B$ can be defined by the matrices
\[
  U_0=\begin{pmatrix}
        0&0&1 \\ 0&0&0 \\ x&y&0
      \end{pmatrix}\ \ 
  U_1=\begin{pmatrix}
        0&0&-y \\ 0&0&x \\ 0&1&0
      \end{pmatrix}\ \ 
  V_0=\begin{pmatrix}
        0&0&0 \\ 0&0&1 \\  y&z&0
      \end{pmatrix}\ \ 
  V_1=\begin{pmatrix}
        0&0&-z \\ 0&0&y \\ -1&0&0
      \end{pmatrix}. 
\]
Note that here $[d]=V_0U_0-U_0V_0$.
\end{example}

\appendix

\section{Graded modules over a non-commutative algebra}
\label{se:appendix}
This section is devoted to the proof of Proposition~\ref{pr:Wmod}, which is a key input in the proof of Theorem~\ref{th:ci_case}. The statement may be viewed as a structure theorem for graded modules over a kind of non-commutative Weyl algebra.

Let $A$ be a commutative ring, $p$ a positive integer, and $A\langle s_1,\dots,s_p,t_1,\dots,t_p\rangle$ the free associative $A$-algebra on indeterminates $\bs s$ and $\bs t$, viewed as a graded $A$-algebra with $|s_i|=1$ and $|t_i|=-1$. We write $[-,-]$ for the graded commutator in a graded ring. Consider the graded $A$-algebra
\[
W\colonequals \frac{A\langle \bs s,\bs t\rangle}{\left([s_i,t_j]=\delta_{ij}\mid 1\le i\le j\le p\right) }\,,
\]
where $\delta_{ij}=1$ if $i=j$ and $0$ if $i\ne j$. One might think of this as the non-commutative, skew-graded, version of the standard Weyl algebra. Below are some basic observations about the product structure of this algebra. To that end, we introduce some notation. Given a sequence $I=(i_1,\dots,i_n)$ of integers with $1\le i_j\le p$ for each $j$ we set
\[
{\bs s}_I \colonequals s_{i_1}\cdots s_{i_n} \quad\text{and} \quad {\bs t}_I \colonequals t_{i_1}\cdots t_{i_n}\,.
\]
In particular, ${\bs s}_{\varnothing}=1={\bs t}_{\varnothing}$. We write $\op I$ for the reverse sequence $(i_n,\dots,i_1)$. We say $I$ is \emph{monotone} to mean that $i_1<\cdots <i_n$; in this case $n\le p$. We write $[p]$ for the monotone sequence $(1,2, \dots, p)$.

\begin{lemma}
\label{le:Walg}
The following statements concerning the algebra $W$ hold.
    \begin{enumerate}[\quad\rm(1)]
        \item Elements of $W$ can be expressed as an $A$-linear combination of monomials
        \[
        s_{i_1}\cdots s_{i_m}t_{j_1}\cdots t_{j_n}\,;
        \]
        and similarly with the order of $\bs s$ and $\bs t$ reversed.
        \item For $1\le i,j\le p$ the elements $s^2_i$ and $[s_i,s_j]$ commute with all elements in $A\langle \bs t\rangle$, and similarly with the roles of $\bs s $ and $\bs t$ switched.
        \item For any subsequence $I\subseteq [p]$ there are equalities
        \[
        \bs t_{\op I} {\bs s_I} = \sum_{J\subseteq I} (-1)^{|J|} {\bs s}_J{\bs t}_{\op J} \quad\text{and}\quad 
                    \bs s_I{\bs t}_{\op I} = \sum_{J\subseteq I} (-1)^{|J|} {\bs t}_{\op J}{\bs s}_J\,, 
        \]
        where  $J$ ranges over the subsequences of $I$ and $|J|$ is its cardinality.
    \end{enumerate}
\end{lemma}

\begin{proof}
Part (1) is immediate from the description of $W$.

(2) To verify that $s_i^2$ commutes with $A\langle \bs t\rangle$ it suffices to verify that $[s_i^2,t_k]=0$ for $1\le k\le p$.  If $i\ne k$, then $s_i\cdot t_k= - t_ks_i$ and hence $s_i^2t_k = t_ks_i^2$. For $i=k$ one has
\[
s_i^2t_i = s_i(1-t_is_i) = s_i - (s_it_i)s_i = s_i - (1-t_is_i)s_i = t_is_i^2\,.
\]
Thus, in either case $[s_i^2,t_k]=0$, as desired. An equally straightforward computation reveals that $[[s_i,s_j],t_k] = [t_k,[s_i,s_j]]$ for all $k$, as desired.

(3) We verify the equality on the left, via an induction on $|I|$; the other equality then follows by symmetry. The basis on the induction is the case $|I|=0$, and then both sides of the desired equality are equal to $1$. Suppose $I=(i_1,\dots,i_n)$ with $n\ge 1$. Setting $\underline{I}=(i_2,\dots,i_{n})$ one gets
\begin{align*}
{\bs t}_{\op I} {\bs s_I} 
    &= {\bs t}_{\op{\underline{I}}} (t_{i_1}s_{i_1}) {\bs s}_{\underline{I}} \\
    &= {\bs t}_{\op{\underline{I}}} (1- s_{i_1}t_{i_1}) {\bs s}_{\underline{I}} \\
    &= {\bs t}_{\op{\underline{I}}}{\bs s}_{\underline{I}} -  ({\bs t}_{\op{\underline{I}}}\cdot s_{i_1})(t_{i_1}\cdot {\bs s}_{\underline{I}})\\
    &={\bs t}_{\op{\underline{I}}}{\bs s}_{\underline{I}} - s_{i_1} ({\bs t}_{\op{\underline{I}}}{\bs s}_{\underline{I}}) t_{i_1}
\end{align*}
where the last equality follows from the fact that $i_1\notin \underline{I}$. Applying the induction hypothesis to $\underline I$ yields the desired equality.
\end{proof}

\begin{definition}
\label{de:exterior-alg}    
Let $E$ be the quotient of $W$ by the left ideal of $W$ generated by elements 
\[
t_i, s_i^2,\text{ and } [s_i,s_j]\quad \text{for $1\le i,j\le p$.}
\]
Thus $E$ is a left $W$-module. Writing $s_i$ also for the residue class of $s_i$ in $W$, it follows from Lemma~\ref{le:Walg} that $E$ can be identified with $A[\bs s]$, the exterior algebra on the $\bs s$, where the $t_i$ act as graded derivations with $t_i\cdot s_j = \delta_{ij}$. 
\end{definition}

The result below may be seen as a structure theorem for graded  modules over the algebra $W$ that are nonzero in only finitely many degrees.

\begin{proposition}
\label{pr:Wmod}
Let $F=\{F_i\}_{i\geqslant 0}$ be a graded left $W$-module with $F_0\ne 0$. Then $F_i\ne 0$ for each $0\le i\le p$. Moreover, if $F_i=0$ for all $i\ge p+1$, then the natural map $W\otimes_A F_0\to F$ induces an isomorphism of left $W$-modules
\[
E \otimes_A F_0 \xrightarrow{\ \cong\ } F\,.
\]
\end{proposition}

\begin{proof}
For degree reasons $t_i\cdot F_0=0$ for all $i$, and hence from Lemma~\ref{le:Walg}(3) we deduce that restricted to $F_0$ one has
\[
{\bs t}_{\op {[p]}} {\bs s}_{[p]} =\mathrm{id}^{F_0}\colon F_0\longrightarrow F_0\,.
\]
Since $F_0\ne 0$ we get that $F_i\ne 0$ for $0\le i\le p$, as claimed. 

From now on suppose $F_i=0$ for $i\ge p+1$. Since $t_i\cdot F_0=0$, restricting the first relation of Lemma~\ref{le:Walg}(3) to $F_1$, we get
\[
0=\id_{F_1}-\sum_{i\in [p]}s_it_i,
\]
hence $F_1\subset A\langle\bs s \rangle F_0$. Restricting successively to $F_2,\dots, F_p$ and arguing by induction, we get $F=A\langle\bs s \rangle F_0$. Using $s_i\cdot F_p=0$ and the second relation of Lemma~\ref{le:Walg}(3) we get similarly the second equality below
\[
A\langle \bs s\rangle F_0 =F = A\langle \bs t\rangle F_p\,;
\]
that is to say, $F$ is generated by $F_0$ even as a module over the subalgebra $A\langle \bs s\rangle$ of $W$, and by $F_p$ as a module over the subalgebra $A\langle \bs t\rangle$ of $W$. Given this and Lemma~\ref{le:Walg}(2) for all $1\le i,j\le p$ one gets
\begin{align*}
t_i^2\cdot F & = t_i^2 \cdot A\langle \bs s\rangle F_0 = A\langle \bs s\rangle t_i^2\cdot F_0=0 \\ 
[t_i,t_j]\cdot F & =   [t_i,t_j]\cdot A\langle \bs s\rangle F_0 = A\langle \bs s\rangle  [t_i,t_j]\cdot F_0=0 \,.
\end{align*}
In the same vein, $s_i^2\cdot F = 0 = [s_i,s_j]\cdot F$. 

It follows from the observations above that the natural map $W\otimes_A F_0\to F$ of graded $W$-modules is onto and that it induces a
map of graded left $W$-modules
\[
\phi\colon E\otimes_A F_0 \longrightarrow F\,.
\]
It remains to verify that this map is one-to-one.

Since $E$ is the exterior algebra over $A$ on the indeterminates $\bs s$, as an $A$-module it is free with basis the monomials $\bs s_I$ where $I$ ranges over the subsequences of $[p]$. Because $\phi$ is a map of graded modules, any element in its kernel of degree $n$ is a finite sum  of the form
$\sum \bs s_I \otimes f_I$ where the $I$ ranges over the subsequences of $[p]$ of size $n$, the  $f_I$ are in $F_0$, and 
\[
\sum \bs s_I f_I = 0 \quad \text{in $F$.}
\]
Fix an ordered sequence $J$ of size $n$. With $J' = [p]\setminus J$ viewed as a sequence, it follows from the relations $s_i^2F=0=[s_i,s_j]F$ that
\[
{\bs s}_{J'}\cdot \sum \bs s_If_I = {\bs s}_{J'}\bs s_Jf_J=\pm {\bs s}_{[p]}f_J\,.
\]
Since ${\bs t}_{\op {[p]}} {\bs s}_{[p]}$ is the identity on $F_0$, we deduce that $f_J=0$. Since $J$ was arbitrary, we conclude that $\sum \bs s_I \otimes f_I=0$, that is to say, $\phi$ is one-to-one, as desired.
\end{proof}


\begin{thebibliography}{10}

\bibitem{Avramov:1998a}
L.~L. Avramov, \emph{{Infinite free resolutions}}, Six lectures on commutative
  algebra (J.~Elias, J.~M. Giral, R.~M. Mir\'o-Roig, and S.~Zarzuela, eds.),
  Progr. in Math., vol. 166, Birkh\"auser Verlag, Basel, 1998.

\bibitem{Avramov:1999a}
\bysame, \emph{{Locally complete intersection homomorphisms and a conjecture of
  Quillen on the vanishing of cotangent homology}}, Ann.\ of Math. \textbf{150}
  (1999), 455--487.

\bibitem{Avramov/Buchweitz/Iyengar/Miller:2010a}
L.~L. Avramov, R.-O. Buchweitz, S.~B. Iyengar, and C.~Miller, \emph{{Homology
  of perfect complexes}}, Adv.\ in Math. \textbf{223} (2010), 1731--1781.

\bibitem{Avramov/Foxby:1991a}
L.~L. Avramov and H.-B. Foxby, \emph{{Homological dimensions of unbounded
  complexes}}, J.~Pure \& Applied Algebra \textbf{71} (1991), no.~2--3,
  129--155.

\bibitem{Avramov/Foxby/Halperin:descent_and_ascent}
L.~L. Avramov, H.-B. Foxby, and S.~Halperin, \emph{Descent and ascent of local
  properties along homomorphisms of finite flat dimension}, J.~Pure \& Applied
  Algebra \textbf{38} (1985), no.~2, 167--185.

\bibitem{Avramov/Foxby/Herzog:1994}
Luchezar~L. Avramov, Hans-Bj{\o}rn Foxby, and Bernd Herzog, \emph{Structure of
  local homomorphisms}, J. Algebra \textbf{164} (1994), no.~1, 124--145.
  \MR{1268330}

\bibitem{Borel/Wallach:2000}
A.~Borel and N.~Wallach, \emph{Continuous cohomology, discrete subgroups, and
  representations of reductive groups}, second ed., Mathematical Surveys and
  Monographs, vol.~67, American Mathematical Society, Providence, RI, 2000.
  \MR{1721403}

\bibitem{Brochard:2017}
Sylvain Brochard, \emph{Proof of de {S}mit's conjecture: a freeness criterion},
  Compos. Math. \textbf{153} (2017), no.~11, 2310--2317. \MR{3692747}

\bibitem{Brochard:2023}
\bysame, \emph{Independent sequences and freeness criteria}, J. Algebra
  \textbf{628} (2023), 486--508. \MR{4575651}

\bibitem{Brochard/Iyengar/Khare:2023a}
Sylvain Brochard, Srikanth~B. Iyengar, and Chandrashekhar~B. Khare, \emph{A
  freeness criterion without patching for modules over local rings}, J. Inst.
  Math. Jussieu \textbf{22} (2023), no.~5, 2117--2129. \MR{4624957}

\bibitem{Brochard/Mezard:2011}
Sylvain Brochard and Ariane M\'ezard, \emph{About de {S}mit's question on
  flatness}, Math. Z. \textbf{267} (2011), no.~1-2, 385--401. \MR{2772256}

\bibitem{Bruns/Herzog:1998a}
W.~Bruns and J.~Herzog, \emph{{Cohen--Macaulay rings}}, Cambridge Studies in
  Advanced Mathematics, vol.~39, Cambridge University Press, 1998, Revised
  edition.

\bibitem{Calegari/Geraghty:2018}
Frank Calegari and David Geraghty, \emph{Modularity lifting beyond the
  {T}aylor-{W}iles method}, Invent. Math. \textbf{211} (2018), no.~1, 297--433.
  \MR{3742760}

\bibitem{StacksProject}
A.~J. de~Jong~et al., \emph{{Stacks project}}, available at\newline
  \url{https://www.math.columbia.edu/~dejong/wordpress/?p=866}.

\bibitem{Foxby:1979a}
Hans-Bj{\o}rn Foxby, \emph{Bounded complexes of flat modules}, J. Pure Appl.
  Algebra \textbf{15} (1979), no.~2, 149--172. \MR{535182}

\bibitem{Iversen:1977}
Birger Iversen, \emph{Amplitude inequalities for complexes}, Ann. Sci. \'Ecole
  Norm. Sup. (4) \textbf{10} (1977), no.~4, 547--558. \MR{568903}

\bibitem{Iyengar/Khare/Manning:2024b}
Srikanth~B. Iyengar, Chandrashekhar~B. Khare, and Jeffrey Manning,
  \emph{Congruence modules and the {W}iles-{L}enstra-{D}iamond numerical
  criterion in higher codimensions}, Invent. Math. \textbf{238} (2024), no.~3,
  769--864. \MR{4824729}

\bibitem{Jorgensen:2010}
Peter J{\o}rgensen, \emph{Amplitude inequalities for differential graded
  modules}, Forum Math. \textbf{22} (2010), no.~5, 941--948. \MR{2719763}

\bibitem{Khare/Thorne:2017}
Chandrashekhar~B. Khare and Jack~A. Thorne, \emph{Potential automorphy and the
  {L}eopoldt conjecture}, Amer. J. Math. \textbf{139} (2017), no.~5,
  1205--1273. \MR{3702498}

\bibitem{Roberts:1987a}
Paul Roberts, \emph{Le th\'eor\`eme d'intersection}, C. R. Acad. Sci. Paris
  S\'er. I Math. \textbf{304} (1987), no.~7, 177--180. \MR{880574}

\bibitem{Taylor/Wiles:1995}
Richard Taylor and Andrew Wiles, \emph{Ring-theoretic properties of certain
  {H}ecke algebras}, Ann. of Math. (2) \textbf{141} (1995), no.~3, 553--572.
  \MR{1333036}

\end{thebibliography}

\newcommand{\noopsort}[1]{}
\providecommand{\bysame}{\leavevmode\hbox to3em{\hrulefill}\thinspace}
\providecommand{\MR}{\relax\ifhmode\unskip\space\fi MR }
\providecommand{\MRhref}[2]{%
  \href{http://www.ams.org/mathscinet-getitem?mr=#1}{#2}
}
\providecommand{\href}[2]{#2}

\end{document}